\documentclass{article}

\usepackage{amsmath,amscd,amsfonts,amssymb,amsthm}
\usepackage{graphicx}
\usepackage{epstopdf}
\usepackage{a4wide}
\usepackage{graphicx}
\usepackage{subfigure}
\usepackage[latin1]{inputenc}
\usepackage{float}
\usepackage{placeins}
\usepackage{multirow}

\newcommand{\keyw}[1]{\par\noindent{\bf Keywords: }#1.}

\begin{document}
\title{An approximation formula for the Katugampola integral}

\author{Ricardo Almeida$^1$\\
{\tt ricardo.almeida@ua.pt}
\and
 Nuno R. O. Bastos$^{1,2}$\\
{\tt nbastos@estv.ipv.pt}}

\date{$^1$Center for Research and Development in Mathematics and Applications (CIDMA)\\
Department of Mathematics, University of Aveiro, 3810--193 Aveiro, Portugal\\
$^2$Department of Mathematics, School of Technology and Management of Viseu \\
Polytechnic Institute of Viseu, 3504--510 Viseu, Portugal}

\maketitle

\def\a{\alpha}
\def\r{\rho}
\def\t{\tau}
\def\DS{\displaystyle}
\def\LI{{I_{a+}^{\a,\r}}}
\def\RI{{I_{b-}^{\a,\r}}}
\keyw{fractional calculus, Katugampola fractional integral, numerical methods}

\begin{abstract}
The objective of this paper is to present an approximation formula for the Katugampola fractional integral, that allows us to solve fractional problems with dependence on this type of fractional operator. The formula only depends on first-order derivatives, and thus we convert the fractional problem into a standard one. With some examples we show the accuracy of the method, and then we present the utility of the method by solving a fractional integral equation.
\end{abstract}

\newtheorem{theorem}{Theorem}[section]
\newtheorem{lemma}[theorem]{Lemma}
\newtheorem{proposition}[theorem]{Proposition}
\newtheorem{corollary}[theorem]{Corollary}
\newtheorem{definition}[theorem]{Definition}

\section{Introduction}

Many real-world phenomena are described by non-integer order systems. In fact, they model several problems since they take into consideration e.g. the memory of the process, friction, flow in heterogenous porous media, viscoelasticity, etc \cite{Benson,Cao,He,Malik,Riewe}. Fractional derivative and fractional integral are generalizations of ordinary calculus, by considering derivatives of arbitrary real or complex order, and a general form for multiple integrals. Although mathematicians have wondered since the very beginning of calculus about these questions, only recently they have proven their usefulness and since then important results have appeared not only in mathematics, but also in physics, applied engineering, biology, etc. One question that is important is what type of fractional operator should be considered, since we have in hand several distinct definitions and the choice dependes on the considered problem. Because of this, we find in the literature several papers dealing with similar subjects, but for different type of fractional operators. So, to overcome this, one solution is to consider general definitions of fractional derivatives and integrals, for which the known ones are simply particular cases. We mention for example the approach using kernels (see \cite{Freed,Malinowska,Odzijewicz1,Odzijewicz2,Srivastava}).

The paper is organized in the following way. In Section \ref{sec:FC} we present the definitions of left and right Katugampola fractional integrals of order $\a>0$. Next, in Section \ref{sec:AF}, we prove the two new results of the paper, Theorems \ref{teo1} and \ref{teo2}. The formula is simple to use, and uses only the function itself and an auxiliary family of functions, where each of them is given by a solution of a Cauchy problem. Finally, in Section \ref{sec:ex}, we present some examples where we compare the exact fractional integral of a test function with some numerical approximations as given in the previous section. To end, we exemplify how we can apply the approximation to solve a fractional integral equation with an initial condition.

\section{Caputo--Katugampola fractional integral}\label{sec:FC}

To start, we review the main concept as presented in \cite{Katugampola1}. It generalizes the Riemann--Liouville and Hadamard fractional integrals by introducing in the definition a new parameter $\r>0$, that allows us to obtain them for special values of $\r$.

\begin{definition} Let $a,b$ be two nonnegative real numbers with $a<b$, $\a,\r$ two positive real numbers, and $x:[a,b]\rightarrow\mathbb{R}$ be an integrable function. The left Katugampola fractional integral is defined as
\begin{equation}\label{def:LI}\LI x(t)=\frac{\r^{1-\a}}{\Gamma(\a)}\int_a^t \t^{\r-1}(t^\r-\t^\r)^{\a-1}x(\t) d\t,\end{equation}
and the right Katugampola fractional integral is defined as
$$\RI x(t)=\frac{\r^{1-\a}}{\Gamma(\a)}\int_t^b \t^{\r-1}(\t^\r-t^\r)^{\a-1}x(\t) d\t.$$
\end{definition}
These notions were motivated from the following relation. When $\a=n$ is an integer, the left Katugampola fractional integral is a generalization of the $n$-fold integrals
$$\int_a^t \t_1^{\r-1}\, d\t_1\,\int_a^{\t_1} \t_2^{\r-1}\, d\t_2\,\ldots \int_a^{\t_{n-1}} \t_n^{\r-1}x(\t_n)\, d\t_n$$
$$=\frac{\r^{1-n}}{(n-1)!}\int_a^t \t^{\r-1}(t^\r-\t^\r)^{n-1}x(\t) d\t.$$
We also notice that, taking $\r=1$, we obtain the left and right Riemann--Liouville fractional integrals, and as $\r\to0^+$, we get the left and right Hadamard fractional integrals (cf. \cite{Kilbas,Podlubny,Samko}). We refer to \cite{Katugampola2}, where a notion of Katugampola fractional derivative is presented, generalizing the Riemann--Liouville and Hadamard fractional derivatives of order $\a\in(0,1)$, and to \cite{Katugampola3} where an existence and uniqueness result is proven for a fractional differential equation involving this new operator.

For example, consider the function
$$x(t)=(t^\r-a^\r)^v, \quad v>-1.$$
Then,
$$\begin{array}{ll}
\DS\LI x(t)& =\DS\frac{\r^{1-\a}}{\Gamma(\a)}\int_a^t \t^{\r-1}(t^\r-\t^\r)^{\a-1}(\t^\r-a^\r)^v d\tau\\
&=\DS\frac{\r^{1-\a}}{\Gamma(\a)}\int_a^t \t^{\r-1}(t^\r-a^\r)^{\a-1}\left[1-\frac{\t^\r-a^\r}{t^\r-a^\r}\right]^{\a-1}(\t^\r-a^\r)^v d\tau.
\end{array}$$
With the change of variables
$$u=\frac{\t^\r-a^\r}{t^\r-a^\r},$$
we arrive to
$$\begin{array}{ll}
\DS\LI x(t) &=\DS \frac{\r^{-\a}}{\Gamma(\alpha)}(t^\r-a^\r)^{\a+v}\int_0^1 (1-u)^{\a-1}u^{v} du\\
&=\DS \frac{\r^{-\a}}{\Gamma(\alpha)}(t^\r-a^\r)^{\a+v}B(\a,v+1),
\end{array}$$
where $B(\cdot,\cdot)$ is the Beta function
$$B(x,y)=\int_0^1 t^{x-1}(1-t)^{y-1}dt, \quad x,y>0.$$
Using the useful property
$$B(x,y)=\frac{\Gamma(x)\Gamma(y)}{\Gamma(x+y)},$$
we get the formula
$$\LI x(t)= \frac{\r^{-\a}\Gamma(v+1)}{\Gamma(\a+v+1)}(t^\r-a^\r)^{\a+v}.$$
In a similar way, if we consider
$$y(t)=(b^\r-t^\r)^v, \quad v>-1,$$
we have the following
$$\RI y(t)= \frac{\r^{-\a}\Gamma(v+1)}{\Gamma(\a+v+1)}(b^\r-t^\r)^{\a+v}.$$

\section{Caputo--Katugampola fractional integral}\label{sec:AF}

In this section, we present the main results of the paper. We prove an approximation formula for the Katugampola fractional integrals, which will allow us later to solve a fractional integral equation, by approximating it by an ordinary differential equation. This idea was motivated by the recent works in \cite{Atan2,Atan1,Pooseh0,Pooseh1}, and has been developed in the recent book \cite{Almeida}, where similar formula are proven for the Riemann--Liouville and Hadamard fractional operators.

\begin{theorem}\label{teo1} Let $N \in \mathbb N$ and  $x:[a,b]\rightarrow\mathbb{R}$ be a function of class $C^1$. For $k \in \{1,\ldots,N\}$, define the quantities
$$A=\DS\frac{\r^{-\a}}{\Gamma(\a+1)}\left[1+\sum_{k=1}^N\frac{\Gamma(k-\a)}{\Gamma(-\a)k!}\right],\quad
B_k=\DS\frac{\r^{1-\a}\Gamma(k-\a)}{\Gamma(\a+1)\Gamma(-\a)(k-1)!},$$
and the function $V_k:[a,b]\to  \mathbb R$ by
$$V_k(t)=\int_a^t\t^{\r-1} (\t^\r-a^\r)^{k-1}x(\t)d\t.$$
Then,
$$\LI x(t)=A(t^\r-a^\r)^\a x(t)-\sum_{k=1}^N B_k(t^\r-a^\r)^{\a-k}V_k(t)+E_N(t),$$
with
$$\lim_{N\to\infty}E_{N}(t)=0.$$
\end{theorem}

\textbf{Proof:} Starting with formula \ref{def:LI}, and integrating by parts choosing
$$\dot{u}(\t)=\t^{\r-1}(t^\r-\t^\r)^{\a-1} \quad \mbox{and} \quad v(\t)=x(\t),$$
we obtain
\begin{equation}\label{eq:1}\LI x(t) =\frac{\r^{-\a}}{\Gamma(\a+1)}(t^\r-a^\r)^{\a}x(a)+\frac{\r^{-\a}}{\Gamma(\a+1)}\int_a^t (t^\r-\t^\r)^\a \dot{x}(\t) d\t.\end{equation}
By the binomial theorem, we have
\begin{equation}\label{eq:2}\begin{array}{ll}
(t^\r-\t^\r)^\a& =\DS \left((t^\r-a^\r)-(\t^\r-a^\r)\right)^\a\\
& =\DS (t^\r-a^\r)^\a\left(1-\frac{\t^\r-a^\r}{t^\r-a^\r}\right)^\a\\
& =\DS (t^\r-a^\r)^\a\sum_{k=0}^\infty \frac{\Gamma(k-\a)}{\Gamma(-\a)k!}\left( \frac{\t^\r-a^\r}{t^\r-a^\r}\right)^k.
\end{array}\end{equation}
Replacing formula \ref{eq:2} into \ref{eq:1}, and if we truncate the sum, we get
$$\LI x(t) =\frac{\r^{-\a}}{\Gamma(\a+1)}(t^\r-a^\r)^{\a}x(a)$$
$$+\frac{\r^{-\a}}{\Gamma(\a+1)}
(t^\r-a^\r)^\a\sum_{k=0}^N \frac{\Gamma(k-\a)}{\Gamma(-\a)k!(t^\r-a^\r)^k}\int_a^t(\t^\r-a^\r)^k\dot{x}(\t) d\t+E_N(t),$$
with
$$E_{N}(t)=\frac{\r^{-\a}}{\Gamma(\a+1)}(t^\r-a^\r)^\a\sum_{k=N+1}^\infty \frac{\Gamma(k-\a)}{\Gamma(-\a)k!}\int_a^t
\left( \frac{\t^\r-a^\r}{t^\r-a^\r}\right)^k\dot{x}(\t) d\t.$$
If we split the sum into $k=0$ and the remaining terms $k=1,\ldots,N$, we deduce the following
$$\LI x(t) =\frac{\r^{-\a}}{\Gamma(\a+1)}(t^\r-a^\r)^{\a}x(t)$$
$$+\frac{\r^{-\a}}{\Gamma(\a+1)}
(t^\r-a^\r)^\a\sum_{k=1}^N \frac{\Gamma(k-\a)}{\Gamma(-\a)k!(t^\r-a^\r)^k}\int_a^t(\t^\r-a^\r)^k\dot{x}(\t) d\t+E_N(t).$$
If we proceed with another integration by parts, choosing this time
$$u(\t)=(\t^\r-a^\r)^k \quad \mbox{and} \quad \dot{v}(\t)=\dot{x}(\t),$$
we obtain
$$\LI x(t) =\frac{\r^{-\a}}{\Gamma(\a+1)}\left[1+\sum_{k=1}^N\frac{\Gamma(k-\a)}{\Gamma(-\a)k!}\right](t^\r-a^\r)^{\a}x(t)$$
$$-\frac{\r^{1-\a}}{\Gamma(\a+1)}\sum_{k=1}^N \frac{\Gamma(k-\a)}{\Gamma(-\a)(k-1)!}(t^\r-a^\r)^{\a-k}
\int_a^t\t^{\r-1} (\t^\r-a^\r)^{k-1}x(\t)d\t+E_N(t),$$
proving the desired formula. It remains to prove that
$$\lim_{n\to\infty}E_{N}(t)=0.$$
Let
$$M=\max_{\tau\in[a,b]}|\dot{x}(\t)|.$$
Since
$$\left( \frac{\t^\r-a^\r}{t^\r-a^\r}\right)^k\leq1,\quad \forall \t\in[a,t],$$
and
$$\sum_{k=N+1}^\infty \left|\frac{\Gamma(k-\a)}{\Gamma(-\a)k!}\right|\leq \sum_{k=N+1}^\infty \frac{\exp(\a^2+\a)}{k^{\a+1}}
\leq\int_N^\infty \frac{\exp(\a^2+\a)}{k^{\a+1}}\, dk= \frac{\exp(\a^2+\a)}{\a N^\a},$$
we get
$$|E_N(t)|\leq \frac{M\r^{-\a}}{\Gamma(\a+1)}(t^\r-a^\r)^\a(t-a) \frac{\exp(\a^2+\a)}{\a N^\a},$$
which converges to zero as $N\to\infty$, ending the proof.

\vspace{1cm}

For the right Katugampola fractional integral, the formula is the following. We omit the proof since it is similar to the proof of Theorem \ref{teo1}.

\begin{theorem}\label{teo2} Let $N \in \mathbb N$ and  $x:[a,b]\rightarrow\mathbb{R}$ be a function of class $C^1$. For $k \in \{1,\ldots,N\}$, define the quantities
$$A=\DS\frac{\r^{-\a}}{\Gamma(\a+1)}\left[1+\sum_{k=1}^N\frac{\Gamma(k-\a)}{\Gamma(-\a)k!}\right], \quad
B_k=\DS\frac{\r^{1-\a}\Gamma(k-\a)}{\Gamma(\a+1)\Gamma(-\a)(k-1)!},$$
and the function $W_k:[a,b]\to  \mathbb R$ by
$$W_k(t)=\int_t^b\t^{\r-1} (b^\r-\t^\r)^{k-1}x(\t)d\t.$$
Then,
$$\RI x(t)=A(b^\r-t^\r)^\a x(t)-\sum_{k=1}^N B_k(b^\r-t^\r)^{\a-k}W_k(t)+E_N(t),$$
with
$$\lim_{N\to\infty}E_{N}(t)=0.$$
\end{theorem}

\section{Examples and applications}\label{sec:ex}

To test the efficiency of the purposed method, consider the test function
$$x(t)=t^{2\r}, \, t\in[0,0.5].$$
The expression of the fractional integral of $x$ is
$${I_{0+}^{\a,\r}} x(t)=\frac{2\r^{-\a}}{\Gamma(\a+3)}t^{\r(\a+2)}, \, t\in[0,0.5].$$

\begin{figure}[ht]
\begin{center}
\subfigure[$\a=0.1$ and $\rho=2.3$]{\label{fig:fig1_ex1}\includegraphics[scale=0.4]{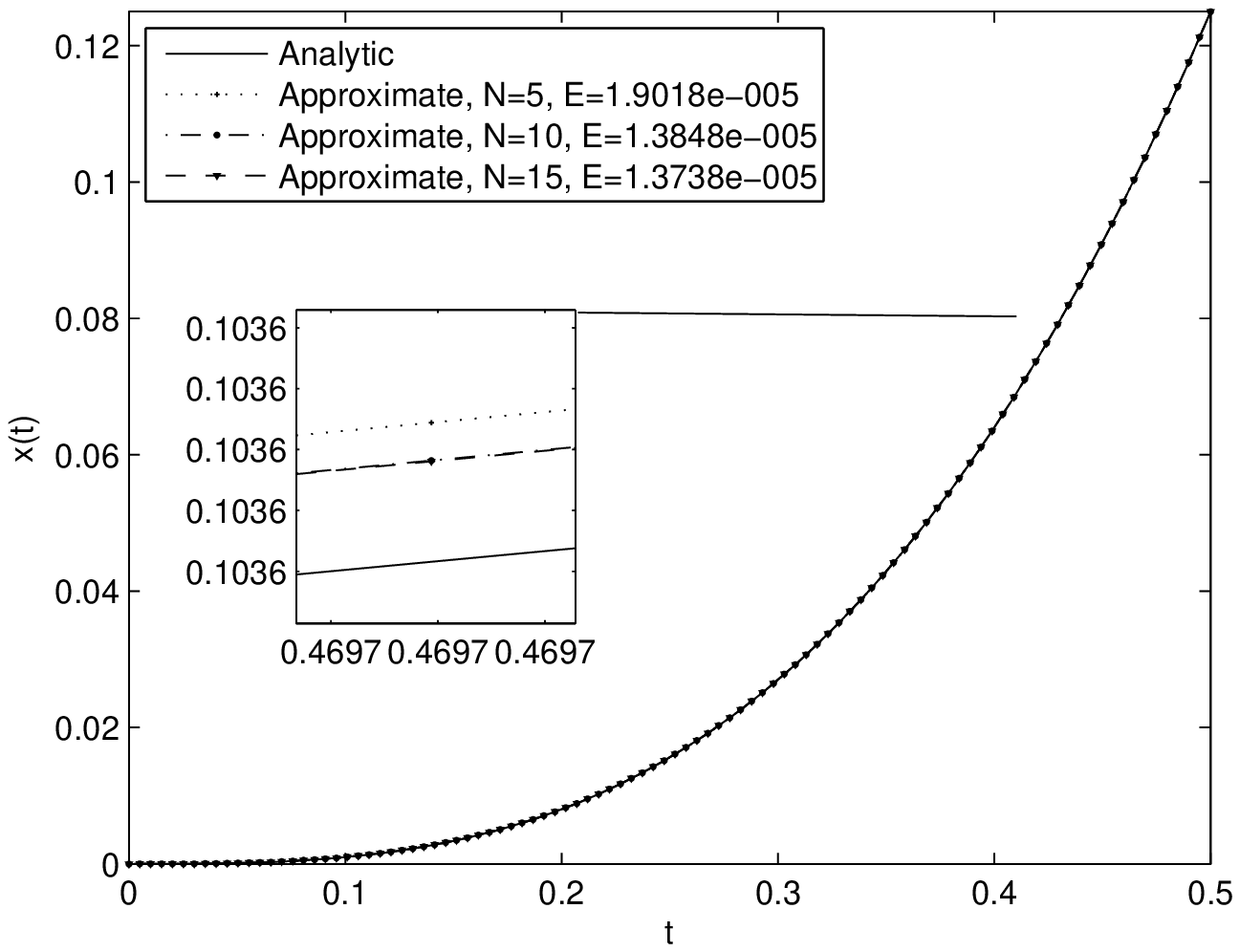}}
\subfigure[$\a=0.9$ and $\rho=0.2$]{\label{fig:fig2_ex1}\includegraphics[scale=0.4]{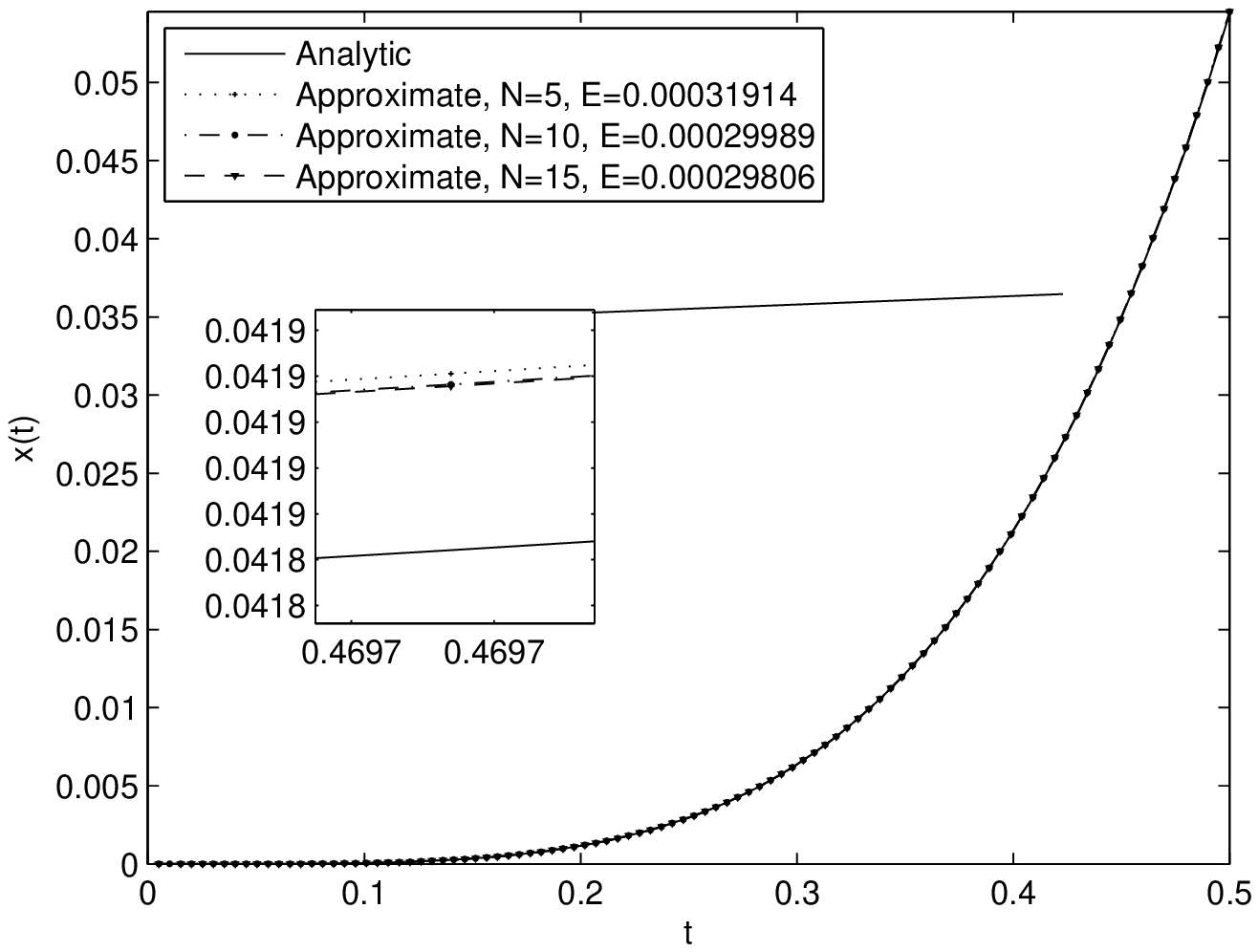}}
\subfigure[$\a=1.5$ and $\rho=0.8$]{\label{fig:fig3_ex1}\includegraphics[scale=0.4]{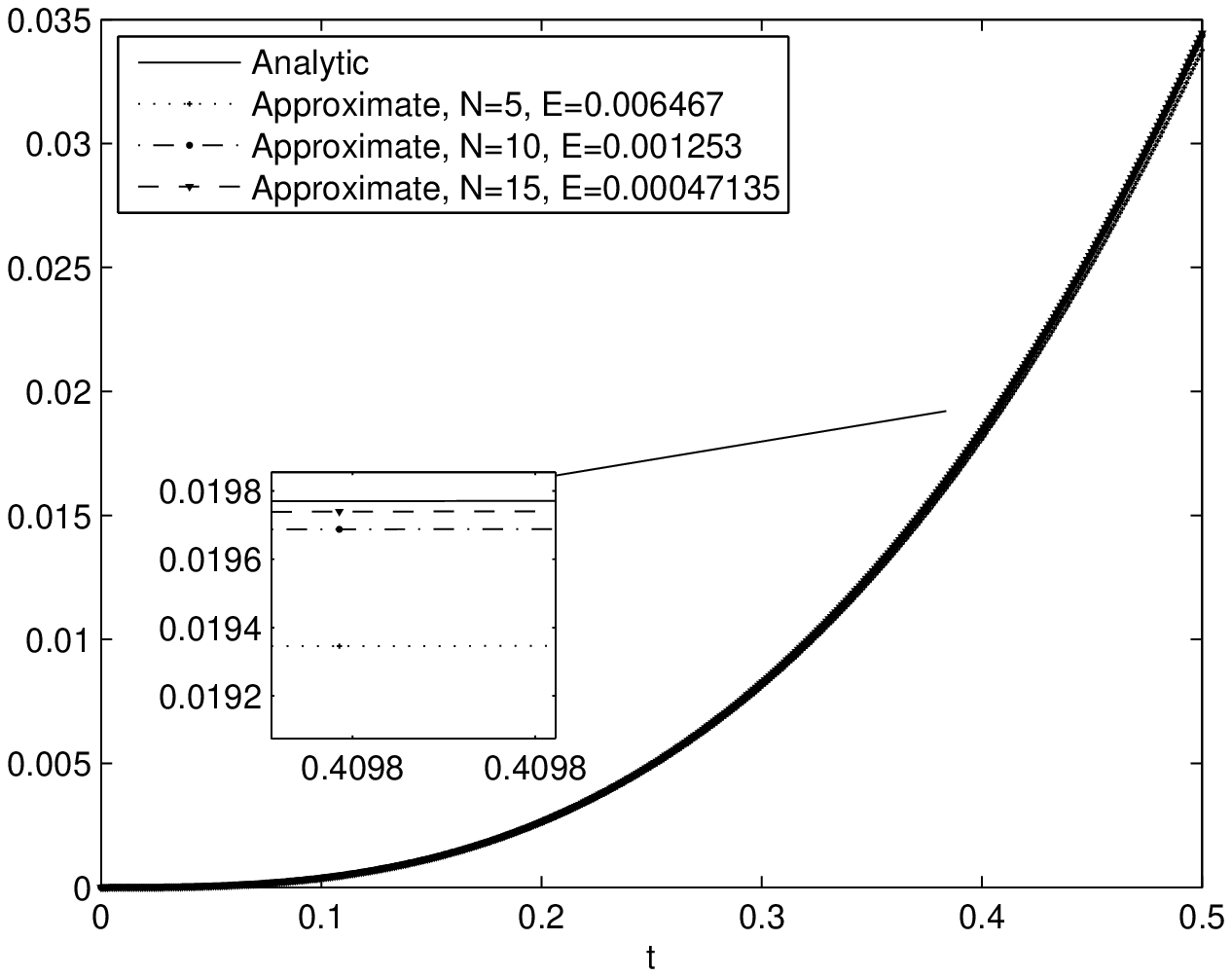}}
\subfigure[$\a=2.4$ and $\rho=0.1$]{\label{fig:fig4_ex1}\includegraphics[scale=0.4]{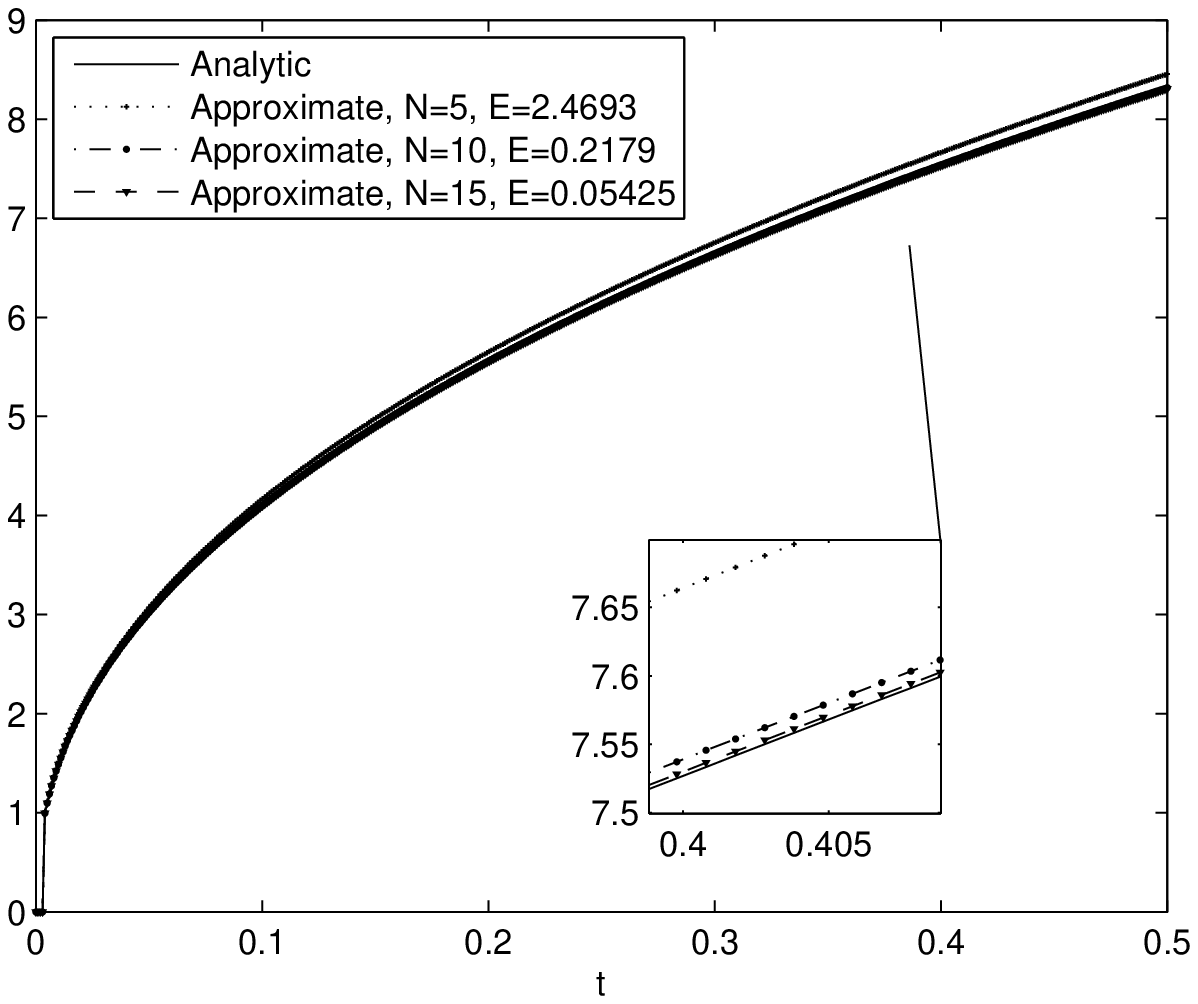}}
\end{center}
\caption{Analytic vs. numerical approximations.}\label{fig11}
\end{figure}

Below, in Figure \ref{fig11}, we present the graphs of the exact expression of the fractional integral of $x$, and some numerical approximations as given by Theorem \ref{teo1} for different values of $\a$ and $\r$.

\begin{figure}[ht]
\begin{center}
\subfigure[$\a=3.5$ and $\rho=1.5$]{\label{fig:fig1_ex2}\includegraphics[scale=0.4]{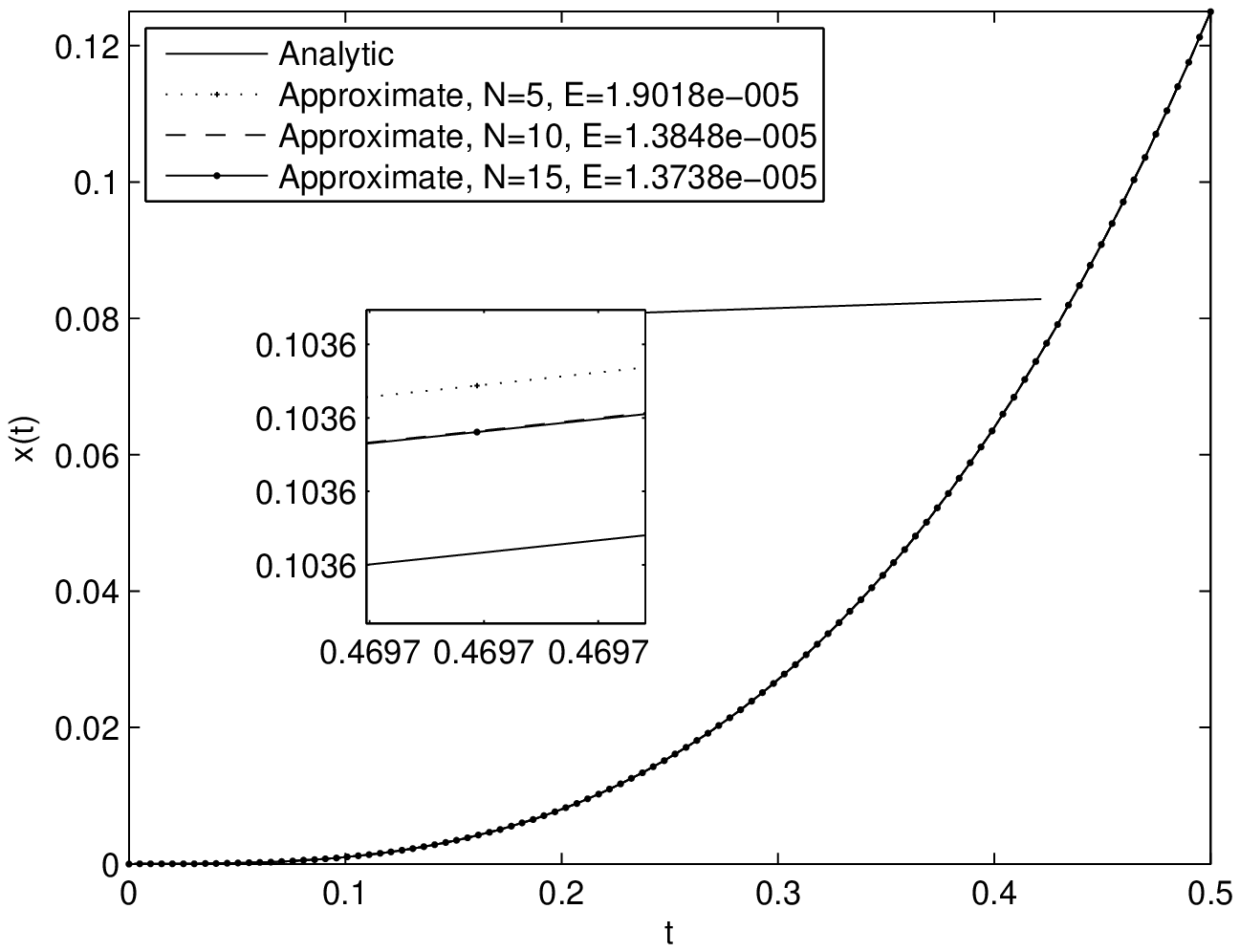}}
\subfigure[$\a=1.8$ and $\rho=2.1$]{\label{fig:fig2_ex2}\includegraphics[scale=0.4]{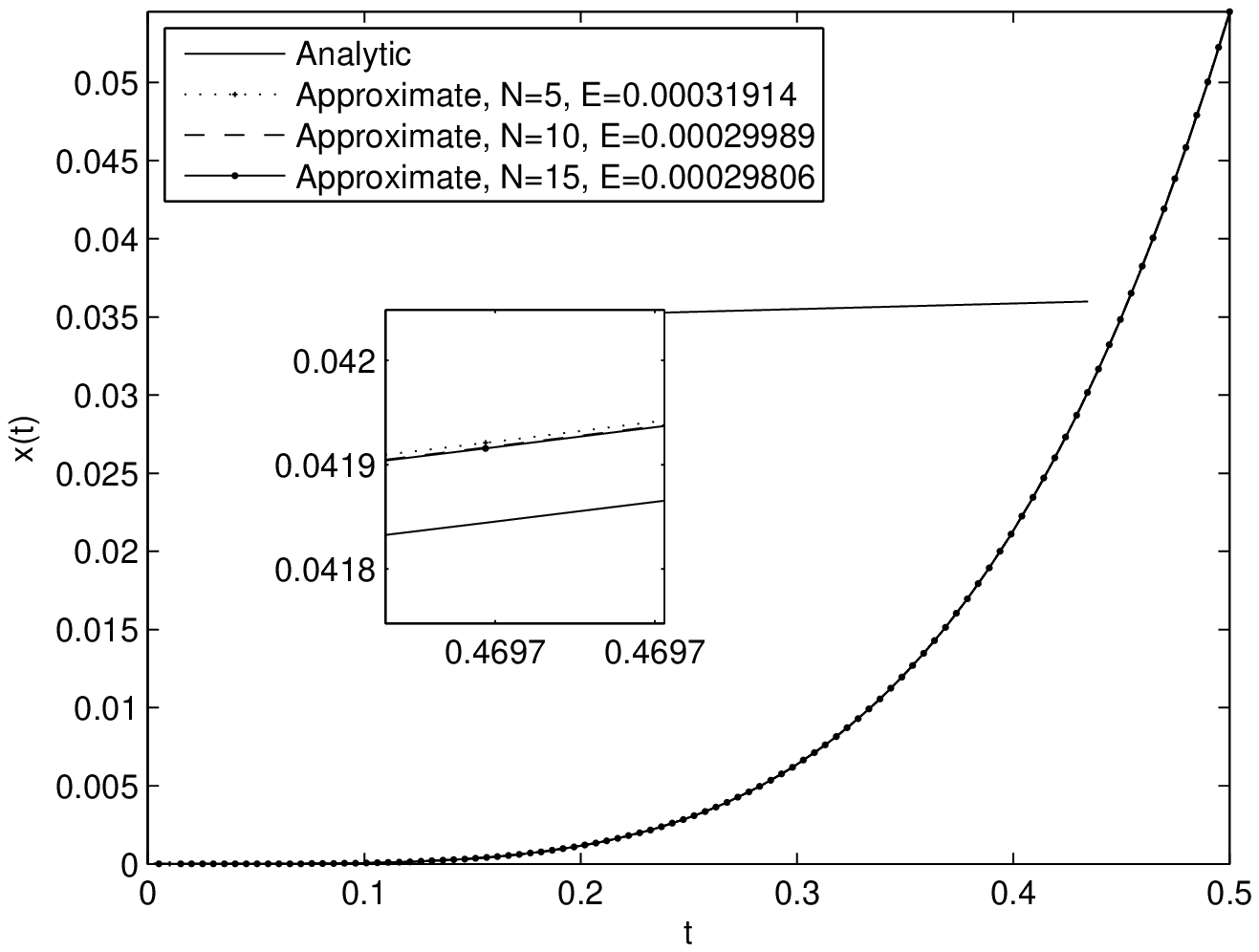}}
\subfigure[$\a=3.3$ and $\rho=1.9$]{\label{fig:fig3_ex2}\includegraphics[scale=0.4]{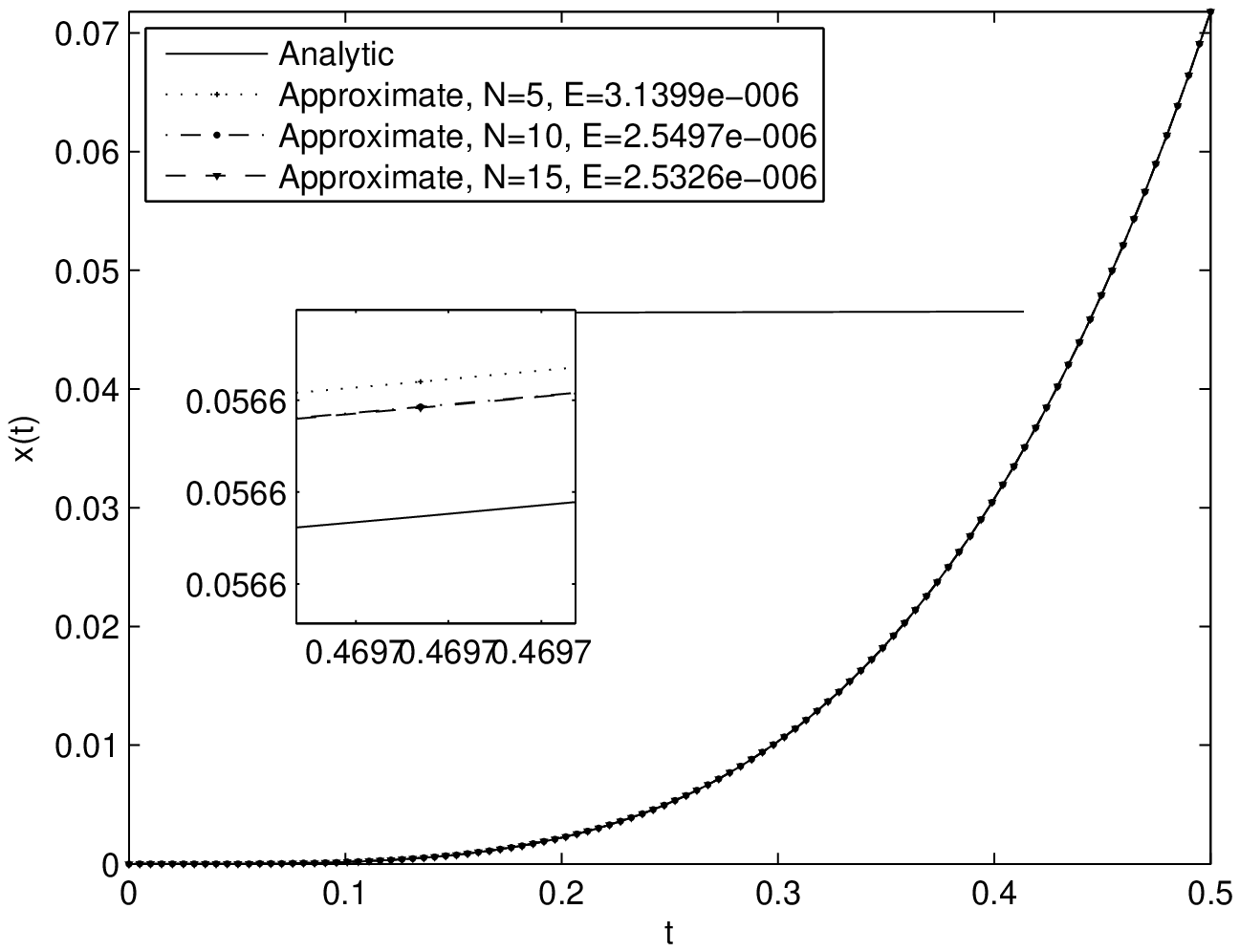}}
\subfigure[$\a=1.4$ and $\rho=1$]{\label{fig:fig4_ex2}\includegraphics[scale=0.4]{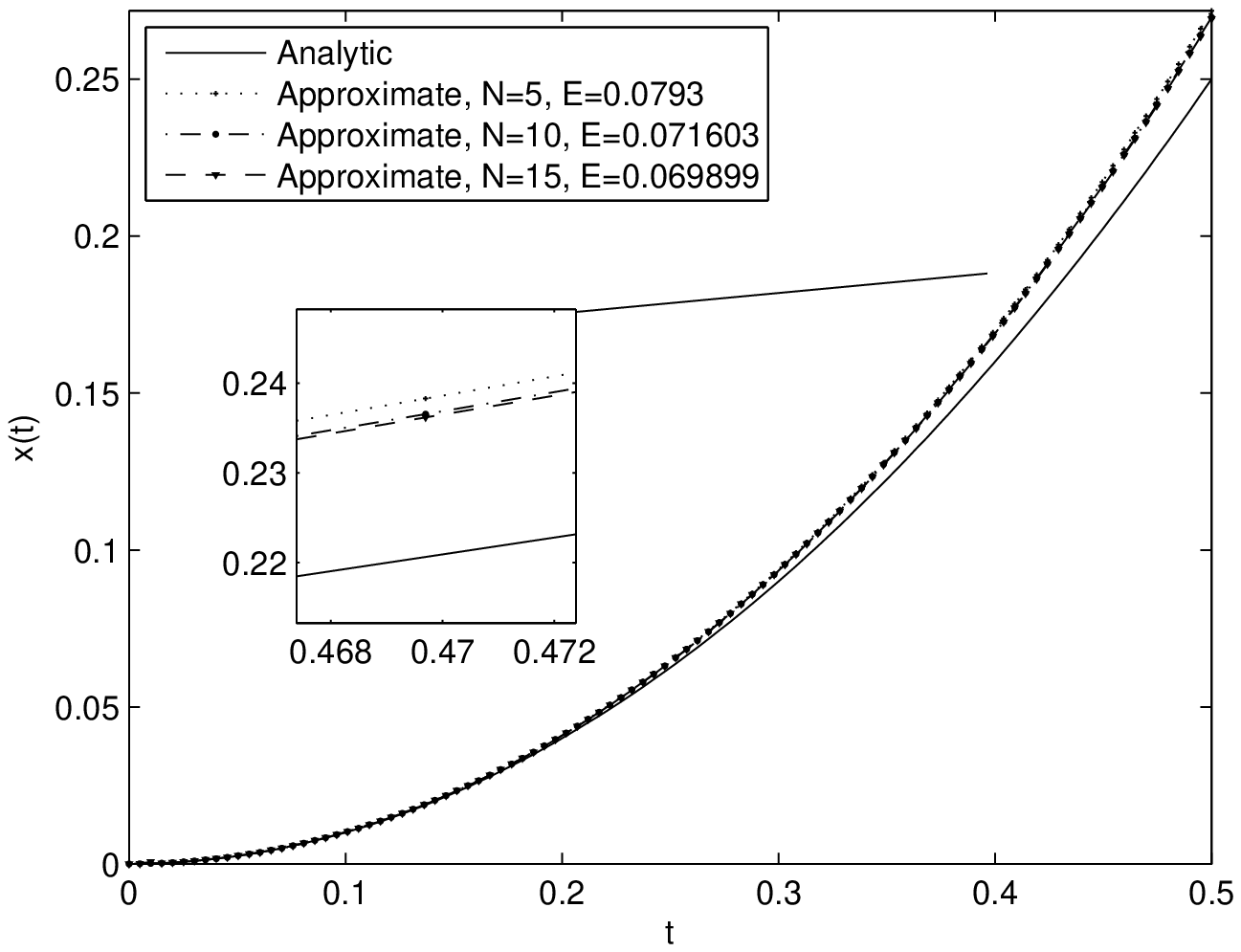}}
\end{center}
\caption{Analytic vs. numerical approximations.}\label{fig12}
\end{figure}

Now we show how the purposed approximation can be useful to solve fractional integral equations with dependence on the Katugampola fractional integral. Consider the system
$$\left\{
\begin{array}{l}
{I_{0+}^{\a,\r}} x(t)+x(t)=t^{2\r}+\frac{2\r^{-\a}}{\Gamma(\a+3)}t^{\r(\a+2)}\\
x(0)=0
\end{array}
\right.$$
The solution for this problem is the function $x(t)=t^{2\r}$. The numerical procedure to solve the problem is the following. Using Theorem  \ref{teo1}, we approximate ${I_{0+}^{\a,\r}} x(t)$ by the sum
$${I_{0+}^{\a,\r}} x(t)\approx A t^{\r\a} x(t)-\sum_{k=1}^N B_k t^{\r(\a-k)}V_k(t),$$
with
$$A=\DS\frac{\r^{-\a}}{\Gamma(\a+1)}\left[1+\sum_{k=1}^N\frac{\Gamma(k-\a)}{\Gamma(-\a)k!}\right],\quad
B_k=\DS\frac{\r^{1-\a}\Gamma(k-\a)}{\Gamma(\a+1)\Gamma(-\a)(k-1)!},$$
and $V_k$ solution of the system
$$\left\{
\begin{array}{l}
\dot{V_k}(t)=t^{\r-1} t^{\r(k-1)}x(t)\\
V_k(0)=0
\end{array}
\right.$$
So, the initial fractional problem is replaced by the Cauchy problem
$$\left\{
\begin{array}{l}
 (A t^{\r\a}+1) x(t)-\sum_{k=1}^N B_k t^{\r(\a-k)}V_k(t)=t^{2\r}+\frac{2\r^{-\a}}{\Gamma(\a+3)}t^{\r(\a+2)}\\
 \dot{V_k}(t)=t^{\r-1} t^{\r(k-1)}x(t), \,\, k=1,\ldots,N\\
x(0)=0\\
V_k(0)=0, \,\, k=1,\ldots,N
\end{array}
\right.$$
For different values of $N$, we obtain different accuracies of the method. Some results are exemplified next, in Figure \ref{fig12}, for different values of $\a$ and $\r$.

\section{Conclusion}
Dealing with fractional operators is in most cases extremely difficult,  and so several numerical methods are purposed to overcome these problems. In our work, we suggest a decomposition formula that depends only on the first-order derivative, and with this tool in hand we can transform the fractional problem into an integer-order one. In Figure \ref{fig11} we considered different values of $\alpha$ and $\rho$, and observe that as $N$ increases, the error of the approximation decreases and the numerical approximations approaches the exact expression of the Caputo--Katugampola fractional integral, converging to it. Next, in Figure \ref{fig12}, we exemplify how it can be useful, by solving a fractional integral equation. For all numerical experiments presented above, we used MatLab to obtain the results.

\section*{Acknowledgments}

This work was supported by Portuguese funds through the CIDMA - Center for Research and Development in Mathematics and Applications,
and the Portuguese Foundation for Science and Technology (FCT-Funda\c{c}\~ao para a Ci\^encia e a Tecnologia), within project UID/MAT/04106/2013.

\end{document}